%% file: multipl-hermite-arxiv.tex
\input neumac
\input multipl-hermite-arxiv.num
\nicknum=0
\def\tr{{\rm tr}}\def\sgn{{\rm sgn\thinspace}}
\nopagenumbers
\let\header=N
\immediate\newwrite\num\immediate\openout\num=multipl-hermite-arxiv.num
\def\leftheadline{\ninepoint\folio\hfill
Multiplier systems}%
\def\transpose#1{\kern1pt{^t\kern-1pt#1}}%

\noindent
{\titelfont Multiplier systems for Hermitian  modular groups}%
\vskip 1.5cm
\centerline{Eberhard Freitag, Universit\"at Heidelberg, Mathematisches Institut}
\centerline{2020}
\vskip1cm
\noindent
\centerline{\vbox{\noindent\hsize=10cm
{\ninepoint{\bf Abstract}%
\smallni
Let $\Gamma_{F,n}$ be the Hermitian modular group of degree $n>1$ in sense of
Hel Braun with respect to an imaginary
quadratic field $F$. Let $r$ be a natural number. There exists a multiplier system of weight $1/r$
(equivalently a Hermitian modular form of weight $k+1/r$, $k$ integral) on some congruence group if
and only if $r=1$ or $2$. This follows from a much more general construction of Deligne [De] combining
it with results of Hill [Hi], Prasad [P] and Prasad-Rapinchuk [PR]. As far as we know, the systems of weight
$1/2$ have not yet been described explicitly. Remarkably
Haowu Wang [Wa] gave an example of a modular form of half integral weight.
Actually he constructs a Borcherds product of weight $23/2$
for a group of type $\O(2,4)$. This group is isogenous to the group $\U(2,2)$ that contains the
Hermitian modular groups of
degree two.
\smallskip
In this paper we want to study such multiplier systems. If one restricts them to the unimodular group
$$\calU=\biggl\{U;\quad\pmatrix{\bar U'^{-1}&0\cr0& U}\in\Gamma\biggr\},$$
one obtains a usual character. Our main result states that the
kernel of this character is a non-congruence subgroup. For sufficiently small $\Gamma$ it coincides
with the group described by Kubota [Ku] in the case $n=2$ and by Bass Milnor Serre [BMS] in the case
$n>2$.
}}}
\vskip 1cm\noindent
{\paragratit Introduction}%
\medskip\noindent
In the paper [FH] we gave a simple proof of a special case of a theorem of Deligne [De1]
that states that the weights of
multiplier systems on subgroups of finite index of the Siegel modular group of degree $n>1$ are integral
or half integral. The same result holds for other modular groups as for the Hermitian modular groups,
the quaternary modular groups and the orthogonal groups $\O(2,n)$, $n\ge 3$. In all these cases multiplier
systems of half integral weight do exist. This follows from results of Deligne [De2], Prasad and Rapunchic
[PR] and Prasad [P]. In the case of the Siegel modular forms these multipliers are obtained
as theta multiplier systems which can be expressed by means of a symplectic Gauss sum. In other cases
no such explicit description is known as far as I know. In this paper we study the multiplier systems
in the Hermitian modular case. Our main result is that the restriction to the unimodular group is
a character whose kernel is a non-congruence subgroup as it has been described in [Ku], [BMS].
\smallskip
Our proof rests on techniques from the paper [FH]. Part of this paper can be generalized from the Siegel
to the Hermitian modular group. In the first part we
generalize results from [FH] to the Hermitian case. Instead of proofs we refer to the corresponding
result in [FH] if the generalization is straightforward.
\smallskip
The main difference of the two cases is that the pair $(R,\gotq)$
for the ring $R$ of integers of an imaginary quadratic field with an distinguished non-zero ideal $\gotq$
admits a non trivial Mennicke symbol (Definition \MenS) whereas every Mennicke symbol in the case $R=\gz$ is trivial.
\smallskip
We fix a natural number $n$ (which later will be 2). We denote by $E=E^{(n)}$ the $n\times n$-unit matrix
and by
$$I=I^{(n)}=\pmatrix{0&-E\cr E&0}$$
the standard alternating matrix.
The unitary group $\U(n,n)$ consists of all
$M\in\GL(2n,\cz)$ with the property $\bar M'IM=I$.
(This is equivalent to $\bar M'HN=H$, where $H$ is the Hermitian form $\imag I$ of signature $(n,n)$).
The special unitary group $\SU(n,n)$ is the subgroup
of elements with determinant one. One has $\SU(1,1)=\SL(2,\rz)$.
\smallskip
From now on we fix an imaginary quadratic field $F=\qz(\sqrt d)$ of discriminant
$d<0$ and denote by
$$\goto=\gz+\gz\omega,\quad \omega={d+\sqrt{d}\over 2},$$
its ring of integers. The Hermitian modular group $\Gamma_{F,n}$ is the subgroup
of $\U(n,n)$ of matrices with entries in $\goto$. Let $\gotq\subset\goto$ be a
non zero ideal. Then
$$\Gamma_{F,n}[\gotq]=\kernel(\Gamma_{F,n}\lo\GL(2n,\goto/\gotq))$$
is the (principal) congruence subgroup of level $\gotq$.
Since the field $F$ is fixed, we can omit the label $F$ and write
$$\Gamma_n:=\Gamma_{F,n}\quad\hbox{and}\quad \Gamma_{n}[\gotq]=\Gamma_{F,n}[\gotq].$$
For  sufficiently small $\gotq$ the group $\Gamma_n[\gotq]$ is contained in $\SL(2n,\goto)$.
Then
the group $\Gamma_1[\gotq]$ is the usual principal congruence subgroup
of the elliptic modular group $\SL(2,\gz)$ of level $\gotq\cap\gz$.
\neupara{Multiplier systems}%
We consider the usual action $MZ=(AZ+B)(CZ+D)^{-1}$ of the unitary group
$\U(n,n)$ on the Hermitian upper half plane
$$\calH_n=\{Z\in\cz^{n\times n};\quad Z=X+\imag Y,\quad X=\bar X',\ Y=\bar Y'>0\ \hbox{(positive definite)}\}.$$
This is an open convex domain in $\cz^{n\times n}$.
The function
$$J(M,Z)=\det(CZ+D)$$
has no zeros on the half plane.
Since the half plane is convex, there exists a continuous choice
$L(M,Z)=\arg J(M,Z)$ of the argument. We normalize it such that
it is the principal value for $Z=\imag E$  where $E$ denotes the
unit matrix.
Recall that the principal value $\Arg(a)$ is defined such that it
is in the interval $(-\pi,\pi]$.
So we have
$$L(M,\imag E)=\Arg(J(M,\imag))\in(-\pi,\pi].$$
We consider
$$w(M,N):={1\over 2\pi}\bigl((L(MN,Z)-L(M,NZ)-L(N,Z)\bigr).$$
Obviously,
$$e^{2\pii w(M,N)}=1.$$ Hence $w(M,N)$ is independent of $Z$ and
$w(M,N)\in \gz$. Usually we will compute $w(M,N)$ by evaluation at $Z=\imag E$.
$$w(M,N)={1\over 2\pi}\bigl(\Arg J(MN,\imag E)-\arg J(M,N(\imag E))-\Arg J(N,\imag E)\bigl),$$
where  $\arg J(M,N(\imag E))$ is obtained through
continuous continuation of the principal value $\Arg J(M,\imag E)$ along a path from $\imag E$ to $N(\imag E)$.
Usually one takes the straight line.
\proclaim
{Remark}
{The function $w:\U(n,n)\times\U(n,n)\to\gz$ is a cocycle
in the following sense:
$$\eqalign{w(M_1M_2,M_3)+w(M_1,M_2)&=w(M_1,M_2M_3)+w(M_2,M_3),\cr
w(E,M)=w(M,E)&=0.\cr}$$
}
Coc%
\finishproclaim
The computation of $w(M,N)$ in degree 1 is easy for the following reason.
From the definition we have
$$2\pi w(M,N)=\Arg((c\alpha+d\gamma)\imag+c\beta+d\gamma)-\arg(c N(\imag)+d)-\Arg(\gamma\imag+\delta)$$
for
$$M=\pmatrix{a&b\cr c&d},\quad N=\pmatrix{\alpha&\beta\cr\gamma&\delta}$$
where $\arg(cN(\imag)+d)$ is obtained from the principal value $\Arg (c\imag +d)$ through
continuous continuation. But $cz+d$ for $z$ in the upper half plane
never crosses the real axis. Hence the
result of the continuation is the principal value too. So all three arguments in the definition
of $w(M,N)$ are the principal values (in degree 1). This makes it easy to compute $w$.
We rely on tables for the values of $w$ which have been derived by Petersson and reproduced by Maass [Ma],
Theorem 16.
\proclaim
{Lemma}
{Let $M={*\;\;\;*\choose m_1\,m_2}$, $S={a\,b\choose c\,d}$ be two real matrices with determinant $1$ and
$(m_1',m_2')$ the second row of the matrix $MS$. Then
$$4w(M,S)=\cases{\sgn c+\sgn m_1-\sgn m_1'-{{\rm sgn}}(m_1cm_1')& if $m_1cm_1'\ne 0$,\cr
-(1-\sgn c)(1-\sgn m_1)& if $cm_1\ne0, m_1'=0$,\cr
(1+\sgn c)(1-\sgn m_2)&if $cm_1'\ne 0,m_1=0$,\cr
(1-\sgn a)(1+\sgn m_1)&if $m_1m_1'\ne 0,c=0$,\cr
(1-\sgn a)(1-\sgn m_2)&if $c=m_1=m_1'=0$.\cr}$$
{\bf Corollary.}
Assume that $m_1cm_1'\ne 0$ and that $m_1m_1'>0$ or $m_1c<0$. Then $w(M,S)=0$.}
MP%
\finishproclaim
We give an example.
\proclaim
{Lemma}
{We have
$$w\biggl(\pmatrix{a&b\cr c&d},\pmatrix{1&x\cr0&1}\biggr)=
w\biggl(\pmatrix{1&x\cr0&1},\pmatrix{a&b\cr c&d}\biggr)=0.$$
}
raTr%
\finishproclaim
\neupara{Some special values of the cocycle}%
We give some examples for values of $w$ in degree $n>1$.
\proclaim
{Lemma}
{One has
$$w\biggl(\pmatrix{E&S\cr0&E},M\biggr)=0.$$}
LTra%
\finishproclaim
The proof is trivial.\qed
\proclaim
{Lemma}
{Let
$$P=\pmatrix{0&1&0&0\cr1&0&0&0\cr 0&0&0&1\cr0&0&1&0}.$$
Set
$$z:=\det(\imag C+D).$$
Then
$$w(P,M)=w(M,P)=\Arg (-z)-\Arg(z)-\Arg(-1).$$}
Pval%
\finishproclaim
{\it Proof.} Compare [FH], Lemma 1.2.\qed
\smallskip
\proclaim
{Definition}
{The {\emph Siegel parabolic group} consists of all elements from $\SU(n,n)$ of the form
$$\pmatrix{A&B\cr0&D}.$$
}
DefSK%
\finishproclaim
There is a  character on the Siegel parabolic group
$$\varepsilon\pmatrix{A&B\cr 0&D}=\det(D).$$
For an element $M$ of the Siegel parabolic group, the expression
$\det(CZ+D)=\det(D)$ is independent of $Z$. Hence
$$L(M,Z)=0\quad\hbox{if}\ \varepsilon(M)=1.$$
An immediate consequence is the following lemma.
\proclaim
{Lemma}
{For two elements $P,Q$ of the Siegel parabolic group we have
$w(P,Q)=0$ if $\varepsilon(P)=1$.}
ParM%
\finishproclaim
The proof is trivial.\qed
\smallskip
We have to consider two embeddings $\iota_1,\iota_2:\SL(2,\rz)\to\SU(2,2)$, namely
$$\iota_1\pmatrix{a&b\cr c&d}=
\pmatrix{a&0&b&0\cr 0&1&0&0\cr c&0&d&0\cr 0&0&0&1},\quad
\iota_2\pmatrix{a&b\cr c&d}=
\pmatrix{1&0&0&0\cr0&a&0&b\cr0&0&1&0\cr0&c&0&d}.$$
\proclaim
{Lemma}
{Let $M$ be in the image of one of the embeddings $\iota_\nu$
and $N$ a Siegel parabolic matrix with $\varepsilon(N)=1$. Then
$w(M,N)=0$.}
KSz%
\finishproclaim
{\it Proof.} Compare [FH], Lemma 3.1.
\proclaim
{Lemma}
{Assume $n=2$. Let
$$M=\pmatrix{E&S\cr0&E},\quad S=\bar S'.$$
Then
$$w(I,M)=\cases{0&if $\tr(S)\ge 0$,\cr -1& else.\cr}$$}
ITra%
\finishproclaim
{\it Proof.}
 Compare [FH], Lemma 1.6.\qed
\proclaim
{Lemma}
{Assume $n=2$. Let
$$M=\pmatrix{E&0\cr S&E},\quad S=\bar S'.$$
Then
$$w(M,I)=\cases{-1&if $\tr(S)> 0$,\cr 0& else.\cr}$$}
TraI%
\finishproclaim
{\it Proof.} Compare [FH], Lemma 1.7.\qed
\neupara{Multipliers}%
\proclaim
{Definition}
{Let $\Gamma\subset\U(n,n)$ be an arbitrary subgroup
and let $r$ be a real number. A system $v(M)$, $M\in\Gamma$, of complex numbers of absolute value $1$
is called a {\emph multiplier system} of weight $r$ if
$$v(MN)\equiv v(M)v(N)\sigma(M,N)$$
where
$$\sigma(M,N)=\sigma_r(M,N):=e^{2\pii r w(M,N)}.$$
}
TriVw%
\finishproclaim
Let now $\Gamma$ be a normal subgroup of finite index of  $\Gamma_n$, $n\ge 2$. Since the congruence subgroup
property holds we know that $\Gamma$ contains a congruence subgroup $\Gamma_n[\gotq]$. It is easy to show that
weights $r$ of multiplier systems are rational [Ch]. Hence a suitable power of $v$ is trivial on some congruence
subgroup. This shows that there exists a natural number $l$ such the all values of $v$ are $l$th roots of unity.
\smallskip
For any $L\in\Gamma_{n}$ we can consider a conjugate
multiplier system on $\Gamma$ that is defined by
$$\tilde v(M)=v(LML^{-1}){\sigma(LML^{-1},L))\over \sigma(L,M)}.$$
It is easy to check that this is a multiplier system and that this defines an action of $\Gamma_{n}$
on the set of all multiplier systems on $\Gamma$.
The quotient of two multiplier systems of the same weight is a homomorphism, as we know into a finite group.
Since the congruence subgroup problem has been solved for the
Hermitian modular group,
we obtain $\tilde v(M)=v(M)$ on some congruence subgroup.
Since the Hermitian modular group is finitely generated, they agree on $\Gamma_n[\gotq]$, $\gotq$ suitable.
\proclaim
{Lemma}
{Let $v$ be a multiplier system on a subgroup $\Gamma\subset\Gamma_n$ of finite index.
In the case $n\ge 2$ there exists an ideal $\gotq\ne 0$ such that $\Gamma[\gotq]\subset\Gamma$
and such that
$$v(M)=v(LML^{-1}){\sigma(LML^{-1},L)\over\sigma(L,M)}\qquad (M\in\Gamma[\gotq])$$
for all $L\in\Gamma_{n}$.
}
EleO%
\finishproclaim
Several times we will replace $\gotq$ by a smaller ideal. We then just say ``for suitable $\gotq$''.
We always assume that suitable $\gotq$ have the property that the only unit of $\goto$ that is congruent to 1 mod $\gotq$
is $1$. Then each Siegel parabolic $M\in\Gamma_n[\gotq]$ has the property $\varepsilon(M)=1$.
We also assume that $\gotq\subset 4\goto$.
\proclaim
{Proposition}
{Let $v$ be a multiplier system on a subgroup $\Gamma\subset\Gamma_2$ of finite index.
For suitable $\gotq$ the group $\Gamma_2[\gotq]$ is contained in $\Gamma$
and for each
matrix $M$ from $\Gamma_2[\gotq]$ of the form
$$M=\pmatrix{E&0\cr*&E}.$$
we have $v(M)=1$.}
ProOe%
\finishproclaim
{\it Proof.} Compare [FH], Proposition 2.4.\qed
\proclaim
{Proposition}
{Let $v$ be a multiplier system on a subgroup $\Gamma\subset\Gamma_2$ of finite index.
For suitable $\gotq$ we have $\Gamma_2[\gotq]\subset\Gamma$ and such the following holds.
Let $U$ be an element from the subgroup that is generated by the matrices
${1\,q\choose 0\,1}$ and ${1\,0\choose q\,1}$ for $q\in\gotq$ and let
$$M=\pmatrix{\bar U'^{-1}&*\cr0&U}.$$
Then $v(M)=1$.}
ProOz%
\finishproclaim
{\it Proof.} Compare [FH], Proposition 2.4.\qed
\neupara{Embedded subgroups}%
We restrict now to the case $n=2$. The case $n>2$ can be derived from this easily.
Besides the embeddings $\iota_1,\iota_2$ we have to consider the embedding
$$\iota:\GL(2,\cz)\lo\U(2,2),\quad
\iota(U)=\pmatrix{\bar U'^{-1}&0\cr0&U}.$$
This gives us an embedding
$\SL(2,\goto)\hookrightarrow \Gamma_2$. We use the notation
$$\SL(2,\goto)[\gotq]=\kernel(\SL(2,\goto)\lo\SL(2,\goto/\gotq)).$$
We have $w(\iota(U),\iota(V))=1$. Hence, for suitable $\gotq$
$$\SL(2,\goto)[\gotq]\lo S^1,\quad U\longmapsto v(\iota(U)),$$
is a homomorphism. We mentioned that the values of $v$ are $l$th roots of unity. Hence the
kernel is a subgroup of finite index in $\SL(2,\goto)[\gotq]$.
\smallskip
Our method depends on some game between the embeddings $\iota_1,\iota_2$ and $\iota$.
We have
$$P\iota_1(M)P^{-1}=\iota_2(M),\quad
P=\pmatrix{0&1&0&0\cr1&0&0&0\cr 0&0&0&1\cr0&0&1&0}.$$
From Lemma \Pval\ follows
$w(\iota_2(M),P)=w(P,\iota_1(M))$. Hence we obtain from Lemma \EleO\
the following result.
\proclaim
{Lemma}
{Let $v$ be a multiplier system on a subgroup $\Gamma\subset\Gamma_2$ of finite index.
For suitable $\gotq$ we have $\Gamma_2[\gotq]\subset\Gamma$
and
$$v(\iota_1( M))=v(\iota_2(M))$$
for $M\in\Gamma_1[\gotq]$.}
iEiZ%
\finishproclaim
For sake of simplicity we write
$$v(M)=v(\iota_1(M))=v(\iota_2(M)).$$
This is a multiplier system  in degree 1.
We have
$$w(M,N)=w(\iota_\nu(M),\iota_\nu(N)),\quad\hbox{for}\quad\nu=1,2.$$
\smallskip
\proclaim
{Lemma}
{Let $v$ be a multiplier system on a subgroup $\Gamma\subset\Gamma_2$ of finite index.
For suitable $\gotq$
the value $v(M)$, $M\in\Gamma_1[\gotq]$, depends only on the second row of~$M$.}
SRd%
\finishproclaim
{\it Proof.} Compare [FH], Lemma 3.2.\qed
\smallskip
\proclaim
{Lemma}
{Let $v$ be a multiplier system on a subgroup $\Gamma\subset\Gamma_2$ of finite index.
For suitable $\gotq$
we have $\Gamma_2[\gotq]\subset\Gamma$ and for any
$$M_1\in \pmatrix{a&b_1\cr c_1&d_1}\in\SL(2,\goto)[\gotq]\quad\hbox{and}
\quad M_2=\pmatrix{a&b_2\cr c_2&d_2}\in\Gamma_1[\gotq]$$
(in particular $a\in\gz$) the relation
$$v\pmatrix{\bar d_1&-\bar c_1&0&0\cr-\bar b_1&a&0&0\cr0&0&a&b_1\cr0&0&c_1&d_1}\cdot
v\pmatrix{a&0&b_2&0\cr0&1&0&0\cr c_2&0&d_2&0\cr0&0&0&1}
= v\pmatrix{1&0&0&0\cr0&a&0&b_1\bar b_1b_2\cr0&0&1&0\cr0&c_1\bar c_1c_2&0&y}$$
holds (where
$y=d_2-b_2c_2(d_1+\bar d_1)+ab_2c_2d_1\bar d_1$).
}
ZweiV%
\finishproclaim
{\it Proof.} Compare [FH], Lemma 3.3.\qed
\smallskip
Now we assume that the multiplier system is of half integral weight. We can restrict it to a subgroup of finite
index of the Siegel modular group. Since this group has the congruence subgroup property, $v$ must agree with
the theta multiplier system on a suitable congruence subgroup. Its restriction to the embedded $\Gamma_1[q]$
is the multiplier system of the classical theta function
$$1+2e^{2\pii z}+2e^{2\pii 4z}+2e^{2\pii 9z}+\cdots,$$
which, for $q\equiv0$ mod 4, is given by the Kronecker symbol  $\bigl({c\over d}\bigr)$(see [Di] for details).
We will need it only for $c\ne 0$ and for
odd $d$.
We collect some properties
(always assuming this condition)
$$\Bigl({c_1c_2\over d}\Bigr)=
\Bigl({c_1\over d}\Bigr)\Bigl({c_2\over d}\Bigr),\quad
\Bigl({c\over d_1d_2}\Bigr)=
\Bigl({c\over d_1}\Bigr)\Bigl({c\over d_2}\Bigr).$$
Assume that $m,n$ are odd coprime numbers such that at least one of them is not negative. Then the usual reciprocity law
$$\Bigl({m\over n}\Bigr)\Bigl({n\over m}\Bigr)=(-1)^{{m-1\over 2}}(-1)^{{n-1\over 2}}$$
holds.
Assume $d>0$ or $c_1c_2>0$. Then
$$\Bigl({c_1\over d}\Bigr)=\Bigl({c_2\over d}\Bigr)\quad\hbox{if}\quad c_1\equiv c_2\mod d.$$
Also the relation
$$\Bigl({c\over d_1}\Bigr)=\Bigl({c\over d_2}\Bigr)\quad\hbox{if}\quad
\cases{d_1\equiv d_2\mod 4c\ \hbox{and}\  c\equiv 0\mod4\ \hbox{or}\cr
d_1\equiv d_2\mod c\ \hbox{and}\ c\equiv 2\mod 4\cr}$$
is valid. Finally we mention
$$\Bigl({c\over -1}\Bigr)=\cases{1&for $c>0$,\cr -1&for $c<0$.\cr}$$
We obtain the existence of a natural number $q\equiv 0$ mod $4$ such that $\gotq\subset(q)$ and such that
$$v(\iota_1 M)=v(\iota_2 M)=\Bigl({c\over d}\Bigr)\quad\hbox{for}\quad M\in\Gamma_1[q].$$
Since $ad\equiv 1$ mod $bc$ and $b\equiv 0$ mod 4, we get from one of the above rules
$$\Bigl({c\over ad}\Bigr)=1\quad\hbox{hence}\quad \Bigl({c\over a}\Bigr)=\Bigl({c\over d}\Bigr).$$
Next we assume that $v$ is trivial on a congruence subgroup inside $\calU$.
So we can get $v(H_2)=1$. From Lemma \ZweiV\ we get the following result.
\smallskip
There exists a natural number $q\equiv 0$ mod $4$ such that for
$$\eqalign{&
a\in\gz,\quad c_1\in\goto,\quad c_2\in\gz,\cr&
a\equiv 1\mod q,\quad c_1\equiv 0\mod q,\quad c_2\equiv 0\mod q,\cr&
a\goto+c_1\goto=\goto,\quad a\gz+c_2\gz=\gz
\cr}$$
the relation
$$\Bigl({c_2\over a}\Bigr)=\Bigl({c_1\bar c_1c_2\over a}\Bigr)$$
holds. This implies
$$\Bigl({c\bar c\over a}\Bigr)=1\quad\hbox{for}\quad a\in 1+q\gz,\ c\in q\goto,\ (a,c)=1.$$
One can apply this relation to $qc$ for an arbitrary $c\in \goto$ to obtain
$$\Bigl({c\bar c\over a}\Bigr)=1\quad\hbox{for}\quad a\in 1+q\gz,\ c\in \goto,\ (a,c)=1.$$
It is known that there are infinitely many primes of the form $p=c\bar c$ [Co]. We choose one such that
$p$ and $q$ are coprime. Then we have
$$\Bigl({p\over a}\Bigr)=\Bigl({a\over p}\Bigr)$$
if $p$ and $a$ are coprime. We choose $\alpha$ such that $\bigl({\alpha\over p}\bigr)=-1$. We can
solve the congruence $1+xq\equiv\alpha$ mod $p$.
Then $a=1+xq$ is coprime to $p$ and we have
$$\Bigl({c\bar c\over a}\Bigr)=
\Bigl({p\over a}\Bigr)=\Bigl({a\over p}\Bigr)=\Bigl({\alpha\over p}\Bigr)=-1.$$
This is a contradiction. This gives the first part of our main results (for $n=2$ and as a consequence
for arbitrary $n$).
\proclaim\
{Theorem}
{Let $\Gamma\subset\Gamma_{n}$ be any
subgroup of finite index of a Hermitian modular group of degree $n\ge 2$. Let $v$ be a multiplier system of
half integral weight. The restriction of $v$ to the subgroup
$$\calU=\biggl\{U;\quad\pmatrix{\bar U'^{-1}&0\cr0& U}\in\Gamma\biggr\}$$
is a usual character. Its kernel is a non-congruence subgroup of finite index.
On a suitable congruence subgroup it agrees with the subgroup constructed in [Ku] in the
case $n=2$ and in [BMS] in the general case.}
ThDz%
\finishproclaim
In the rest of the paper we will give the proof of the second part. (The case $n=2$ is enough.)
\neupara{Mennicke symbol}%
We recall the notion of a Mennicke symbol.
Let $R$ be a Dedekind domain and $\gotq\subset R$ a non-zero ideal.
We introduce the set
$$\calC(R,\gotq):=\bigl\{\>(a,b)\in R\times R;\quad Ra+Rb=R,\
a\equiv 1\mod\gotq,\ b\equiv 0\mod \gotq\>\bigr\}.$$
Every pair $(a,b)$ is the first column of a matrix $\bigl({a\,c\atop b\,d}\bigr)\in\SL(2,R)$.
Multiplying $M$ with   a matrix of the type $\bigl({1\,x\atop 0\,1}\bigr)\in\SL(2,R)$ one can achieve
$c\equiv 0$ mod $\gotq$ and $d\equiv 1$ mod $\gotq$.
\proclaim
{Definition}
{A Mennicke symbol mod $\gotq$ is a map
$$\calC(R,\gotq)\lo G,\quad (a,b)\loma\Bigl[{b\atop a}\Bigr],$$
into some group $G$
such that the following properties hold.
\smallni
{\rm MS1} It is invariant under the transformations $(a,b)\mapsto (a+xb,b)$ and $(a,b)\mapsto (a,b+qay)$
for integral $x,y$.
\smallni
{\rm MS2} It satisfies the rule
$$\Bigl[{b_1b_2\atop a}\Bigl]=\Bigl[{b_1\atop a}\Bigl]\Bigl[{b_2\atop a}\Bigl].$$}
MenS%
\finishproclaim
In our context, the group $G$ will be the group of complex numbers of absolute value one.
Mennicke symbols have been classified in [BSM] for  Dedekind domains of arithmetic type.
If $R$ is the ring of algebraic integers in a number field that is not totally imaginary, then
the Mennicke symbols are trivial. In the case of a totally imaginary field they can be described explicitly
by means of power residue symbols.
\smallskip
The main result of this section is the following theorem.
\proclaim
{Theorem}
{Let $v$ be a multiplier system of half integral weight on a subgroup of finite index
of a Hermitian modular group of degree two. Then there exists a non-zero  ideal $\gotq\subset\goto$
with the following properties.
\smallni
{\rm 1)} $\Gamma_2[\gotq]\subset\Gamma$.
\vskip1mm\noindent
{\rm 2)} There exists a Mennicke symbol\/ $[\cdot]$ for $(\goto,\gotq)$
such that for all $U\in\SL(2,\goto)[\gotq]$ one has
$$\Bigl[{c\atop a}\Bigr]=v\pmatrix{\bar U'^{-1}&0\cr 0&U},\quad U=\pmatrix{a&b\cr c&d}.$$}
PisM%
\finishproclaim
{\it Proof.} We mention that the kernel of $v$ on $\SL(2,\goto)[\gotq]$ agrees for suitable $\gotq$ with a
non congruence subgroup constructed by Kubota [Ku].
The proof of the theorem is given during the rest of this section.
\smallskip
We have to consider also the embeddings $\iota_1,\iota_2:\Gamma_1[q]\lo\Gamma_2[q]$.
As in the Hermitian case we have
$v(i_1(M))=v(i_2(M))$ an this depends only on the second row of $M\in\Gamma_1[q]$. Hence we
can define
$$\Bigl\{{c\atop d}\Bigr\}=v\pmatrix{a&0&b&0\cr 0&1&0&0\cr c&0&d&0\cr 0&0&0&1}^{-1}=
v\pmatrix{1&0&0&0\cr0&a&0&b\cr0&0&1&0\cr0&c&0&d}^{-1}.$$
The elements of $\calC(\goto,\gotq)$ are the second columns of the matrices in $\GL(2,\goto)[\gotq]$.
Hence
$$\Bigl[{b\atop d}\Bigr]=\biggl(\iota\pmatrix{a&b\cr c&d}\biggr)$$
is well-defined on $\calC(\goto,\gotq)$. We claim that this symbol satisfies MS1.
We notice that $w$ is trivial on the image of $\iota$. Hence $v$ is a character on this group.
The invariance under $(a,b)\mapsto (a,b+qay)$ follows from the equation
$$\pmatrix{a&b\cr c&d}\pmatrix{1&qy\cr 0&1}=\pmatrix{a&b+qay\cr*&*}.$$
To prove the invariance under $(a,b)\mapsto (a+xb,b)$, we consider
$$\pmatrix{1&0\cr-x&1}\pmatrix{a&b\cr c&d}\pmatrix{1&0\cr x&1}=\pmatrix{a+xb&b\cr*&*}.$$
Due to Lemma \EleO\ we can assume that $v(\iota(M)))$ is invariant under conjugation
with $\iota{1\,0\choose x\,1}$. This proves MS1.
\smallskip
We would like to have also MS2. To get a result in this direction, we make use of
$$v(\iota_\nu(M^{-1}))=v(\iota_\nu(M))^{-1},\quad\nu=1,2.$$
This is true since in degree 1 one has $w(M,M^{-1})=0$. (This is a general rule for $c\ne0$
and also for $c=0$ and $a>0$. But in our case $c=0$ implies $a=1$ since we assume $q>2$.)
From the analogue of Lemma \ZweiV\ we get  the general rule
(compare Lemma 13.3 in [BMS].)
$$\Bigl[{c_1\atop a}\Bigr]\Bigl\{{c_2\atop a}\Bigr\}=\Bigl\{{c_1\bar c_1c_2\atop a}\Bigr\}.$$
We insert $c_2=1-a$.
\proclaim
{Lemma}
{We have
$$\Bigl\{{1-a\atop a}\Bigr\}=1$$
for $a\equiv 1\mod q$.}
emaA%
\finishproclaim
{\it Proof.} We use
$$\pmatrix{1&1\cr0&1}\pmatrix{1&0\cr a-1&1}\pmatrix{1&-1\cr0&1}=\pmatrix{2-a&a-1\cr 1-a&a}.\eqno\square$$
Now we obtain
$$\Bigl[{c\atop a}\Bigr]=\Bigl\{{c\bar c(1-a)\atop a}\Bigr\}.$$
Before we continue, we mention that $\{\}$ is not  a Mennicke symbol. It does not satisfy MS1.
Nevertheless it is closely related to $[\cdot]$.
\proclaim
{Lemma}
{We have
$$\Bigl\{{c\atop d}\Bigr\}=\Bigl\{{c\atop d+yc}\Bigr\}$$
and
$$\Bigl\{{c+xqd\atop d}\Bigr\}=\Bigl\{{c\atop d}\Bigr\}e^{2\pii rs}\quad
\hbox{where}\quad s=w\biggl(\pmatrix{*&*\cr c&d},\pmatrix{1&0\cr qx&1}\biggr).$$
}
GnM%
\finishproclaim
{\it Proof.}
The first relation
can be derived from
$$\pmatrix{1&-y\cr 0&1}\pmatrix{*&*\cr c&d}\pmatrix{1&y\cr0&1}=\pmatrix{*&*\cr c&d+cy}.$$
To derive the second one we consider the relation
$$\pmatrix{*&*\cr c&d}\pmatrix{1&0\cr qx&1}=\pmatrix{*&*\cr c+dxq&d}.$$
It shows
$$\Bigl\{{c+dxq\atop b}\Bigr\}=\Big\{{c\atop d}\Bigr\}\;e^{2\pii rs}.$$
The $w$-value $s$ is usually not zero.\qed
\smallskip
But from the corollary of the table of Maass in the introduction we get
$$w\biggl(\pmatrix{*&*\cr c\bar c&a},\pmatrix{1&0\cr-qc\bar c&1}\biggl)=0.$$
Using this, we get
$$\Bigl\{{c\bar c(1-a)\atop a}\Bigr\}=\Bigl\{{c\bar c\atop a}\Bigr\}.$$
So we obtain
$$\Bigl[{c\atop a}\Bigr]=\Bigl\{{c\bar c\atop a}\Bigr\}$$
and moreover
$$\Bigl[{c_1c_2\atop a}\Bigr]=\Bigl\{{c_1\bar c_1c_2\bar c_2\atop a}\Bigr\}=
\Bigl[{c_1\atop a}\Bigr]\Bigl\{{c_2\bar c_2\atop a}\Bigr\}
=\Bigl[{c_1\atop a}\Bigr]\Bigl[{c_2\atop a}\Bigr].$$
This is part of the condition MS2. (We assume up to now $a\in\gz$).
\proclaim
{Lemma}
{Let $(a,b)$ be two elements of $\goto$ such that $(a,b)=\goto$. Then there exists $x\in\goto$
such that $a+xb$ is not divisible by any natural number $>1$.}
abX%
\finishproclaim
{\it Proof.} In the case $b=0$ we can take $x=0$. Hence we can assume
that $b\ne 0$.
\smallskip
We write an element $a\in\goto$ in the form
$$a=\dot a+\ddot a\omega,\quad\omega={d+\sqrt d\over 2}.$$
We will use
$$\omega^2=-N(\omega)+d\omega.$$
From a solution $ax+by=1$ we derive that the 4 integers
$$\dot a,\ \dot b,\ \ddot a N(\omega),\ \ddot bN(\omega)$$
are coprime. We have to find $x\in\goto$ such that
$$\dot a+\dot x\dot b-\ddot x\ddot bN(\omega),\quad
\ddot a+\dot x\ddot b+\ddot x(\dot b+\ddot b d)$$
are coprime. We consider the greatest common divisors
$$g=\ggT(\ddot a,\dot b,\ddot b)$$
By Dirichlet's prime number theorem we can find $y\in\goto$ such that
$$\ddot a+\dot y\ddot b+\ddot y(\dot b+\ddot b d)=gp$$
where $p$ is a prime number. There are infinitely many choices for $p$. Hence we can
get that $p$ is coprime to $\ddot b^2(d^2-d)/4+\dot b^2+\dot b\ddot b d$.
(This expression equals $N(\dot b+\ddot b\omega)$ which is positive.)
Now we set
$$\dot x=\dot y+t(\dot b+\ddot bd),\quad\ddot x=\ddot y-t\ddot b\qquad(t\in\gz).$$
Then we have still
$$\ddot a+\dot x\ddot b+\ddot x(\dot b+\ddot b d)=gp$$
but
$$\dot a+\dot x\dot b-\ddot x\ddot bN(\omega)=
\dot a+\dot y\dot b-\ddot y\ddot bN(\omega)+t(\dot b^2+\dot b\ddot b d+\ddot b^2N(\omega)).$$
Now we consider the greatest common divisor
$$g'=(\dot a+\dot y\dot b-\ddot y\ddot bN(\omega),\dot b^2+\dot b\ddot b d+\ddot b^2N(\omega)).$$
We can choose $t$ such that
$$\dot a+\dot x\dot b-\ddot x\ddot bN(\omega)=g'p'$$
where $p'$ is a prime. We can choose $p'$ coprime to $gp$. Our goal was to get
$gp$ and $g'p'$ coprime. This means that $g,g'$ are coprime. But
$$\ggT(g,g')=\ggT(\ddot a,\dot b,\ddot b,\dot a+\dot y\dot b-
\ddot y\ddot bN(\omega),\dot b^2+\dot b\ddot b d+\ddot b^2N(\omega))=\ggT(\dot a,\dot b,\ddot a,\ddot b)=1.$$
This proves Lemma \abX.\qed
\smallskip
Now we come to the proof of MS2. We can write it in the form
$$\Bigl[{q^2b_1b_2\atop a}\Bigl]=\Bigl[{qb_1\atop a}\Bigl]\Bigl[{qb_2\atop a}\Bigl]$$
where $a,b_1,b_2$ are in $\goto$ such that $a\equiv 1$ mod $q$ and such that
$(a,b_1)=(a,b_2)=\goto$.
This formula is invariant under the replacement
$b_1\mapsto b_1+xa$. By Lemma \abX\ 
we can assume that
$b_1$ is not divisible by any natural number. We also want to make an replacement for $b_2$.
For this we consider the ray class of the principal ideal $(b_2)$ mod the ideal $(a)$.
(Recall that two ideals $\gotb_1$, $\gotb_2$ are in the same ray class mod an ideal $\gota$
if there exist $\beta_1\equiv\beta_2\equiv1$ mod $\gota$ such that $\beta_1\gotb_1=\beta_2\gotb_2$.)
Our product formula does not change if one replaces $b_2$ by $\beta b_2$ for $\beta\equiv 1$ mod $(a)$.
Hence we may replace $(b_2)$ by any other $(b_2')$ in the same ray class. In each ray class there
are infinitely many primes. Hence we can assume that $b_2$ is  coprime
to $N(b_1)$. Now we make the replacement $b_2+xaN(b_1)$. Again we make use of Lemma \abX\
to reduce
to the case where $b_2$ is not divisible by any natural number, and, in addition, is coprime to $Nb_1$.
We claim that
then also $b_1b_2$ is not divisible by any natural number. We argue  indirectly and assume that there is
prime number $p$ that divides $b_1b_2$. Then $p$ splits in $\goto$ into two prime ideals,
$(p)=\bar{\gotp}\gotp$. Since $p$ does not divide $b_1$ and $b_2$, we can assume $\bar{\gotp}\vert b_1$ and
$\gotp\vert b_2$. But then $\gotp\vert\bar b_1$ which contradicts to the fact that $\bar b_1$
and $b_2$ are coprime.
\smallskip
Now we make the stronger assumption $a\equiv 1$ mod $q^2$. Then $\ddot a\equiv0$ mod $q^2$.
We can make the replacement
$a\mapsto a+xq^2b_1b_2$ without changing the product formula. This means that we replace
$\ddot a\loma q^2(\ddot a/q^2+\ddot y)$ where $y=xb_1b_2$. Since $b=b_1b_2$ is not divisible by any
natural number, $\ddot y$ runs through all rational integers if $x$ runs through all
integers in $\goto$. This shows that can transform $\ddot a$ to $0$. So
we can assume $a\in\gz$. But in this case the product formula has been proved. Now we replace $q$ by $q^2$
to obtain Theorem \PisM.\qed
\vfill\eject\noindent
{\paragratit References}
\bigskip
\item{[BMS]} Bass, H. Milnor, J. Serre, J.P.: {\it Solution of the congruence subgroup problem
for $\SL_n\; (n\ge 3)$ and $\Sp_{2n}\; (n\ge 2)$,}
Publications math\`ematiques l'I.H.\'E.S., tome {\bf 33}, p. 59--137 (1967)
\medskip
\item{[Co]} Cox,\ D.A.: {\it Primes of the Form $x^2+ny^2$,}
Wiley Series in Pure and Applied Mathematics (2013)
\medskip
\item{[De1]} Deligne, P.: {\it Extensions centrales non r\'esiduellement finies de groupes
arithmetiques,} C. R. Acad. Sci. Paris {\bf 287}, p. 203--208 (1978)
\medskip
\item{[De2]} Deligne, P.: {\it Extensions Centrales  de groupes Alg\'ebriques
Simplememt Connexes et Cohomology Galoisienne,}
Publ. Math. I.H.E.S. {\bf84}, 35--89 (1996)
\medskip
\item{[Di]} Dickson, L.E.: {\it Introduction to the Theory of Numbers,}
Dover Publications, New York, Dover (1957).
\medskip
\item{[FH]} Freitag, E., Hauffe-Waschb\"usch, A.:
{\it Multiplier systems for Siegel modular groups,}
\hbox{arXiv: 2009.06455} [math.NT] (2020)
\medskip
\item{[Hi]} Hill, R.: {\it Fractional weights and non-congruence subgroups,} Automorphic
Forms and Representations of algebraic groups over local fields,
Saito,\ H., Takahashi,\ T. (ed.) Surikenkoukyuroku series {\bf 1338}, 71--80 (2003)
\medskip
\item{[Ku]} Kubota, T.: {\it Ein arithmetischer Satz \"uber eine Matrizengruppe,}
J.\ reine angew.\ Math., {\bf 222}, 55--77 (1965)
\medskip
\item{[Ma]} Maass, H.: {\it Lectures on Modular Functions of One Complex Variable,} Notes by Sunder
Lal, Tata Institute Of Fundamental Research, Bombay, Revised 1983  (1964)
\medskip
\item{[Me]} Mennicke, J.: {\it Zur Theorie der Siegelschen Modulgruppe,}
Math.\ Annalen {\bf 159}, 115--129 (1965)
\item{[P]} Prasad, P.:
{Deligne's topological central extension is universal},
Advances in Mathematics {\bf 181}, 160--164 (2004)
\medskip
\item{[PR]} Prasad, G., Rapinchuk, A.: {\it Computation of the metaplectic
kernel,} Publications mathe\'ematiques de l'I.H.\'E.S. tome {\bf 84},
p.91--187 (1996)
\bye

%% file: neumac.tex
%
\output={\if N\header\headline={\hfill}\fi
\plainoutput\global\let\header=Y}
\magnification\magstep1
\tolerance = 500
\hsize=14.4true cm
\vsize=22.5true cm
\parindent=6true mm\overfullrule=2pt
\newcount\kapnum \kapnum=0
\newcount\parnum \parnum=0
\newcount\procnum \procnum=0
\newcount\nicknum \nicknum=1
\font\ninett=cmtt9

\font\ninebf=cmbx9

\font\sixbf=cmbx6
\font\ninesl=cmsl9

\font\nineit=cmti9

\font\ninerm=cmr9

\font\sixrm=cmr6
\font\ninei=cmmi9
\font\eighti=cmmi8
\font\sixi=cmmi6
\skewchar\ninei='177 \skewchar\eighti='177 \skewchar\sixi='177
\font\ninesy=cmsy9
\font\eightsy=cmsy8
\font\sixsy=cmsy6
\skewchar\ninesy='60 \skewchar\eightsy='60 \skewchar\sixsy='60
\font\titelfont=cmr10 scaled 1440
\font\paragratit=cmbx10 scaled 1200

\font\name=cmcsc10
\font\emph=cmbxti10

\font\tenmsbm=msbm10
\font\sevenmsbm=msbm7
%

%
\font\got=eufm10

\font\teneufm=eufm10
\font\seveneufm=eufm7
\font\fiveeufm=eufm5
\newfam\eufmfam
\textfont\eufmfam=\teneufm
\scriptfont\eufmfam=\seveneufm
\scriptscriptfont\eufmfam=\fiveeufm

\font\tenmsam=msam10
\font\sevenmsam=msam7
\font\fivemsam=msam5
\newfam\msamfam
\textfont\msamfam=\tenmsam
\scriptfont\msamfam=\sevenmsam
\scriptscriptfont\msamfam=\fivemsam
\font\tenmsbm=msbm10
\font\sevenmsbm=msbm7
\font\fivemsbm=msbm5
\newfam\msbmfam
\textfont\msbmfam=\tenmsbm
\scriptfont\msbmfam=\sevenmsbm
\scriptscriptfont\msbmfam=\fivemsbm
\def\Bbb#1{{\fam\msbmfam\relax#1}}
\def\cz{{\kern0.4pt\Bbb C\kern0.7pt}
}
\def\ez{{\kern0.4pt\Bbb E\kern0.7pt}
}
\def\fz{{\kern0.4pt\Bbb F\kern0.3pt}}
\def\gz{{\kern0.4pt\Bbb Z\kern0.7pt}}
\def\hz{{\kern0.4pt\Bbb H\kern0.7pt}
}
\def\kz{{\kern0.4pt\Bbb K\kern0.7pt}
}
\def\nz{{\kern0.4pt\Bbb N\kern0.7pt}
}
\def\oz{{\kern0.4pt\Bbb O\kern0.7pt}
}
\def\rz{{\kern0.4pt\Bbb R\kern0.7pt}
}
\def\sz{{\kern0.4pt\Bbb S\kern0.7pt}
}
\def\pz{{\kern0.4pt\Bbb P\kern0.7pt}
}
\def\qz{{\kern0.4pt\Bbb Q\kern0.7pt}
}
\newskip\ttglue
\def\ninepoint{\def\rm{\fam0\ninerm}%
  \textfont0=\ninerm \scriptfont0=\sixrm \scriptscriptfont0=\fiverm
  \textfont1=\ninei \scriptfont1=\sixi \scriptscriptfont1=\fivei
  \textfont2=\ninesy \scriptfont2=\sixsy \scriptscriptfont2=\fivesy
  \textfont3=\tenex \scriptfont3=\tenex \scriptscriptfont3=\tenex
  \def\it{\fam\itfam\nineit}%
  \textfont\itfam=\nineit
  \def\sl{\fam\slfam\ninesl}%
  \textfont\slfam=\ninesl
  \def\bf{\fam\bffam\ninebf}%
  \textfont\bffam=\ninebf \scriptfont\bffam=\sixbf
   \scriptscriptfont\bffam=\fivebf
  \def\tt{\fam\ttfam\ninett}%
  \textfont\ttfam=\ninett
  \tt \ttglue=.5em plus.25em minus.15em
  \normalbaselineskip=11pt
  \font\name=cmcsc9
  \let\sc=\sevenrm
  \let\big=\ninebig
  \setbox\strutbox=\hbox{\vrule height8pt depth3pt width0pt}%
  \normalbaselines\rm
  \def\sl{\it}}

\headline={\ifodd\pageno\rightheadline\else\leftheadline\fi}
\def\rightheadline{\ninepoint Paragraphen"uberschrift\hfill\folio}
\def\leftheadline{\ninepoint\folio\hfill Chapter"uberschrift}
\let\header=Y
\def\titel#1{\need 9cm \vskip 2truecm
\parnum=0\global\advance \kapnum by 1
{\baselineskip=16pt\lineskip=16pt\rightskip0pt
plus4em\spaceskip.3333em\xspaceskip.5em\pretolerance=10000\noindent
\titelfont Chapter \uppercase\expandafter{\romannumeral\kapnum}.
#1\vskip2true cm}\def\leftheadline{\ninepoint
\folio\hfill Chapter \uppercase\expandafter{\romannumeral\kapnum}.
#1}\let\header=N
}
\def\Titel#1{\need 9cm \vskip 2truecm
\global\advance \kapnum by 1
{\baselineskip=16pt\lineskip=16pt\rightskip0pt
plus4em\spaceskip.3333em\xspaceskip.5em\pretolerance=10000\noindent
\titelfont\uppercase\expandafter{\romannumeral\kapnum}.
#1\vskip2true cm}\def\leftheadline{\ninepoint
\folio\hfill\uppercase\expandafter{\romannumeral\kapnum}.
#1}\let\header=N
}
\def\need#1cm {\par\dimen0=\pagetotal\ifdim\dimen0<\vsize
\global\advance\dimen0by#1 true cm
\ifdim\dimen0>\vsize\vfil\eject\noindent\fi\fi}
\def\neupara#1{\par\penalty-2000
\procnum=0\global\advance\parnum by 1
\vskip1cm\noindent{\paragratit \the\parnum. #1}%
\def\rightheadline{\ninepoint\S\the\parnum.\ #1\hfill \folio}%
\vskip 8mm\noindent}
\def\Proclaim #1 #2\finishproclaim {\bigbreak\noindent
{\bf#1\unskip{}. }{\it#2}\medbreak\noindent}
%
\gdef\proclaim #1 #2 #3\finishproclaim {\bigbreak\noindent%
\global\advance\procnum by 1
{%
{\relax\ifodd \nicknum
\hbox to 0pt{\vrule depth 0pt height0pt width\hsize
   \quad \ninett#3\hss}\else {}\fi}%
\bf\the\parnum.\the\procnum\ #1\unskip{}. }
{\it#2}
\immediate\write\num{\string\def
 \expandafter\string\csname#3\endcsname
 {\the\parnum.\the\procnum}}
\medbreak\noindent}
\newcount\stunde \newcount\minute \newcount\hilfsvar
\def\uhrzeit{
    \stunde=\the\time \divide \stunde by 60
    \minute=\the\time
    \hilfsvar=\stunde \multiply \hilfsvar by 60
    \advance \minute by -\hilfsvar
    \ifnum\the\stunde<10
    \ifnum\the\minute<10
    0\the\stunde:0\the\minute~Uhr
    \else
    0\the\stunde:\the\minute~Uhr
    \fi
    \else
    \ifnum\the\minute<10
    \the\stunde:0\the\minute~Uhr
    \else
    \the\stunde:\the\minute~Uhr
    \fi
    \fi
    }

\def\calC{{\cal C}} 
 
 \def\calH{{\cal H}}

\def\calU{{\cal U}}

\def\gota{\hbox{\got a}} 
\def\gotb{\hbox{\got b}}

\def\gotp{\hbox{\got p}} 
\def\gotq{\hbox{\got q}} 
\def\goto{\hbox{\got o}}

\def\Arg{\mathop{\rm Arg}\nolimits}

\def\ggT{\mathop{\rm ggT}\nolimits}
\def\GL{\mathop{\rm GL}\nolimits}

\def\kernel{\mathop{\rm kernel}\nolimits}

\def\mod{\mathop{\rm mod}\nolimits}
\def\O{{\rm O}}
\def\U{{\rm U}}

\def\SL{\mathop{\rm SL}\nolimits}

\def\SU{\mathop{\rm SU}\nolimits}
\def\Sp{\mathop{\rm Sp}\nolimits}

\def\boxit#1{
  \vbox{\hrule\hbox{\vrule\kern6pt
  \vbox{\kern8pt#1\kern8pt}\kern6pt\vrule}\hrule}}
\def\Boxit#1{
  \vbox{\hrule\hbox{\vrule\kern2pt
  \vbox{\kern2pt#1\kern2pt}\kern2pt\vrule}\hrule}}

\def\smallni{\smallskip\noindent }

\def\lo{\longrightarrow}

\def\loma{\longmapsto}
\def\imag{{\rm i}}
\def\pii{\pi {\rm i}}

\def\square{\hbox{\hbox to 0pt{$\sqcup$\hss}\hbox{$\sqcap$}}}
\def\qed{\ifmmode\square\else{\unskip\nobreak\hfil
\penalty50\hskip3em\null\nobreak\hfil\square
\parfillskip=0pt\finalhyphendemerits=0\endgraf}\fi}
\def\pn{\the\parnum.\the\procnum}
\def\downmapsto{{\buildrel
        {\vbox{\hbox{\hskip.2pt$\scriptstyle-$}}}
        \over{\raise7pt\vbox{\vskip-4pt\hbox{$\textstyle\downarrow$}}}}}